\input amstex
\documentstyle{amsppt}
\input bull-ppt
\keyedby{bull237/lic}
\define\a{\operatorname{Aut}}
\define\te{\operatorname{trace}}
\define\hm{\operatorname{Hom}}
\define\ot{\operatorname{tr}}
\define\br{\bold R}
\define\bc{\bold C}
\define\bz{\bold Z}
\define\bq{\bold Q}
\define\wt{\widetilde}
\define\ov{\overline}

\topmatter
\cvol{26}
\cvolyear{1992}
\cmonth{Jan}
\cyear{1992}
\cvolno{1}
\cpgs{87-112}
\title $\Lambda$\<-Trees and Their Applications \endtitle
\author John W. Morgan\endauthor
\shortauthor{J. W. Morgan}
\shorttitle{$\Lambda$\<-Trees and Their Applications}
\address Department of Mathematics, Columbia University,
New York, New York 10027 \endaddress
\ml jm\@math.columbia.edu\endml
\date March 13, 1991\enddate
\thanks This paper was given as a Progress in Mathematics 
Lecture at the August
8--11, 1990 meeting in Columbus, Ohio\endthanks
\endtopmatter

\document
To most mathematicians and computer scientists the word 
``tree'' conjures
up, in addition to the usual image, the image of a 
connected graph with
no circuits. We shall deal with various aspects and 
generalizations of
these mathematical trees. (As Peter Shalen has pointed 
out, there will
be leaves and foliations in this discussion, but they do 
not belong to
the trees!) In the last few years various types of trees 
have been the
subject of much investigation. But this activity has not 
been exposed
much to the wider mathematical community. To me the 
subject is very
appealing for it mixes very na\"\i ve\ geometric 
considerations with 
the very sophisticated geometric and algebraic 
structures. In fact,
part of the drama of the subject is guessing what type of 
techniques
will be appropriate for a given investigation: Will it be 
direct and
simple notions related to schematic drawings of trees or 
will it be
notions from the deepest parts of algebraic group theory, 
ergodic
theory, or commutative algebra which must be brought to 
bear? Part
of the beauty of the subject is that the na\"\i ve\ tree 
considerations
have an impact on these more sophisticated topics. In 
addition, trees
form a bridge between these disparate subjects.

Before taking up the more exotic notions of trees, let us 
begin with
the graph-theoretic notion of a tree. A graph has 
vertices and edges
with each edge having two endpoints each of which is a 
vertex. A 
graph is a simplicial tree if it contains no loops. In 
\S1 we shall
discuss simplicial trees and their automorphism groups. 
This study is
closely related to combinatorial group theory. The later 
sections of 
this article shall focus on generalizations of the notion 
of a simplicial
tree. A simplicial tree is properly understood to be a 
$\bz$\<-tree. For
each ordered abelian group $\Lambda$, there is an 
analogously defined
object,  called a $\Lambda$\<-tree. A very important 
special case is
when $\Lambda=\br$. When $\Lambda$ is not discrete, the 
automorphisms
group of a $\Lambda$\<-tree is no longer combinatorial in 
nature. One
finds mixing (i.e., ergodic) phenomena occurring, and the 
study of the
automorphisms is much richer and less well understood. We 
shall outline
this more general theory and draw parallels and contrasts 
with the
simplicial case.

While the study of trees and their automorphism groups is 
appealing
{\it per se}, interest in them has been mainly generated 
by considerations
outside the subject. There are several basic properties 
of trees that
account for these connections. As we go along we shall 
explain these
notions and their connections in more detail but let me 
begin with an
overview.

\subheading{\RM1.\quad One-dimensionality of trees} Trees 
are clearly of dimension
one. This basic property is reflected in relations of 
trees to both
algebraic and geometric  objects.

\ch \RM{(a)} Algebraic aspects\endch The resurgence of 
the study of trees 
began with Serre's book \cite{20} on $SL_2$ and trees. In 
this book,
Serre showed how there is naturally a tree associated to 
$SL_2(K)$
where $K$ is a field with a given discrete valuation. The 
group
$SL_2(K)$ acts on this tree. More generally, if $G(K)$ is 
a semisimple
algebraic group of real rank 1 (like $SO(n,1)$\<) and if 
$K$ is a field
with a valuation with formally real residue field, then 
there is a
tree on which $G(K)$ acts. The reason for this can be 
summarized by 
saying that the local Bruhat-Tits building \cite2 for a 
rank 1 group is 
a tree. The fact that the group is rank 1 is reflected in 
the fact that
its local building is one-dimensional. (One does not need 
to understand
to entire Bruhat-Tits machinery to appreciate this 
connection. In fact,
one can view this case as the simplest, yet 
representative, case of the
Bruhat-Tits theory.)

\ch\RM{(b)} Codimension-\RM1 dual objects\endch Suppose 
that $M$ is a manifold
whose fundamental group acts on a tree. Since the tree is 
1-dimensional
there is a dual object in the manifold which is 
codimension-1. This 
object turns out to be a codimension-1 lamination with a 
transverse
measure. These dynamic objects  are closed related to 
isometry groups
of trees. In fact, dynamic results for these laminations 
can be used to
study group actions on trees. They also give rise to many 
interesting
examples of such actions.

\subheading{\RM2.\quad Negative curvature of trees} 
Suppose that we have a
simply connected Riemannian manifold with distance 
function $d$ of
strictly negative curvature. In $M$ the following 
geometric property
holds: If $A_1,A_2,B_1,B_2$ are points with $d(A_1,A_2)$ 
and $d(B_1,
B_2)$ of reasonable size but $d(A_1,B_1)$ extremely 
large, the geodesics
$\gamma_1$ and $\gamma_2$ joining $A_1$ to $B_1$ and 
$A_2$ to $B_2$ are
close together over most of their length. (See Figure 1.) 
The estimate
on how close will depend on an upper bound for the 
curvature and will
go to 0 as this bound goes to $-\infty$. In fact, Gromov 
\cite8 has
defined a class of negatively curved metric spaces in 
terms of this
4-point property. From this point of view $\br$\<-trees 
are the most
negative curved of all spaces, having curvature 
$-\infty$, since in a
tree the geodesics $\gamma_1$ and $\gamma_2$ coincide 
over most of their
lengths. 

We state this fact in another way. Suppose that we have a 
sequence on
strictly negatively curved, simply connected manifolds 
with curvature
upper bounds going $-\infty$, then a subsequence of these 
manifolds
converges in a geometric sense to an $\br$\<-tree. If the 
manifolds in
question are the universal coverings of manifolds with a 
given fundamental
group $G$, the actions of $G$ of the manifolds in the 
subsequence
converge to an action of $G$ on the limit tree. The space 
of negatively
curved manifolds with fundamental group $G$ is completed 
by adding ideal
points at infinity which are represented 
by actions of $G$ on $\br$\<-trees.
This 

\fighere{8 pc}\caption{Figure 1}

\noindent is of particular importance for the manifolds 
of constant negative
curvature---the hyperbolic manifolds. This of course 
brings us 
full circle since the group of automorphisms of 
hyperbolic $n$\<-space is 
the real rank-1 group $SO(n,1)$.

The paper is organized in the following manner. The first 
section is
devoted to discussing the now classical case of 
simplicial trees and
their relationship to combinatorial group theory and to 
$SL_2$ over a 
field with a discrete valuation. Here we follow \cite{20} 
closely. This
material is used to motivate all other cases.

Section 2 begins with the definition of a 
$\Lambda$\<-tree for a general
ordered abelian group $\Lambda$. We show how if $K$ is a 
field with a
nondiscrete valuation with value group $\Lambda$, then 
there is an action
of $SL_2(K)$ on a $\Lambda$\<-tree. More generally, for 
any semisimple
algebraic group of real rank 1 there is a similar result. 
We explain
the case of $SO(n,1)$ in some detail. The last part of 
the section
discusses some of the basics of the actions of groups on 
$\Lambda$ trees.
We introduce the hyperbolic length of a single isometry 
and define and
study some of the basic properties of the space of all 
projective classes
of nontrivial actions. At the end of the section we 
discuss the concept
of base change. 

Section 3 discusses the application of this material to 
compactify the
space of conjugacy classes of representations of a given 
finitely
presented group into a rank-1 group such as $SL_2$ or 
$SO(n,1)$. We
approach this both algebraically and geometrically. The 
main idea is
that a sequence of representations of a fixed finitely 
generated group
$G$ into $SO(n,1)$ has a subsequence which converges 
modulo conjugation
either to a representation into $SO(n,1)$ or to an action 
of $G$ on an
$\br$\<-tree. In the case when the representations are 
converging to an 
action on a tree, this result can be interpreted as 
saying that as one
rescales hyperbolic space by factors going to zero, the 
actions of $G$
on hyperbolic space coming from the sequence of 
representations of $G$
into $SO(n,1)$ converge to an action of $G$ on an 
$\br$\<-tree. This can
then be used to compactify the space of conjugacy classes 
of representations
of $G$ in $SO(n,1)$ with the ideal points being 
represented as actions of
$G$ on $\br$\<-trees.

Section 4 we consider the relationship between 
$\br$\<-trees and 
codimension-1 measured laminations. We define a 
codimension-1 
measured lamination in a manifold. We show how `most' 
measured
laminations in a manifold $M$ give rise to actions of 
$\pi_1(M)$
on $\br$\<-trees. We use these to give examples of 
actions of groups
on $\br$\<-trees, for example, actions of surface groups. 
We also
show how any action of $\pi_1(M)$ on an $\br$\<-tree can 
be dominated
by a measured lamination.

In the last section we give applications of the results 
from the previous
sections to study the space of hyperbolic structures 
(manifolds of
constant negative curvature) of dimension $n$ with a 
given fundamental
group $G$. We have a compactification of this space where 
the ideal points
at infinity are certain types of actions of $G$ on 
$\br$\<-trees. We
derive some consequences---conditions under which the 
space of hyperbolic
structures is compact. The theme running through this 
section is that to
make the best use of the results of the previous sections 
one needs to 
understand the combinatorial group-theoretic information 
that can be
derived from an action of a group on an $\br$\<-tree. For 
simplicial
actions this is completely understood, as we indicate in 
\S1. For more
general actions there are some partial results (coming 
mainly by using
measured laminations), but no general picture. We end the 
article with 
some representative questions about groups acting on 
$\br$\<-trees and
give the current state of knowledge on these questions.

 For other introductions to the subject on group actions 
on $\Lambda$\<-trees
see \cite{21, 12, 3, or 1}. 

\heading 1. Simplicial trees,
 combinatorial group theory, and $SL_2$\endheading

In this discussion of the `classical case' of simplicial 
trees, we are
following closely Serre's treatment in {\it Trees\/} 
\cite{20}. By an
abstract 1-complex we mean a set $V$, the set of 
vertices, and a set $E$,
the set of oriented edges. The set $E$ has a free 
involution $\tau\:
e\mapsto \ov e$, called {\it reversing the orientation}. 
There is also 
a map 
$$
\partial\: E\to V\times V,\qquad \partial e=(i(e),t(e)),
$$
which associates to each oriented edge its initial and 
terminal vertex
satisfying $i(\ov e)=t(e)$. Closely related to an 
abstract 1-complex is
its {\it geometric realization}, which is a topological 
space. It is
obtained by forming the disjoint union
$$
V\coprod\coprod_{e\in E} I_e,
$$
where each $I_e$ is a copy of the closed unit interval, 
and taking the
quotient space under the following relations:
\roster
\item "(i)" $0\in I_e$ is identified with $i(e)\in V$. 
\item "(ii)" $1\in I_e$ is identified with $t(e)\in V$. 
\item "(iii)" $s\in I_e$ is identified with $1-s$ in 
$I_{\ov e}$.
\endroster
The {\it order\/} of a vertex $v$ is the number of 
oriented edges
$e$ for which $i(e)=v$. 

A 1-complex is finite if it has finitely many edges and 
vertices. In
this case its geometric realization is a compact space. 
Figure 2 gives
a typical example of a finite 1-complex.

A nonempty, connected 1-complex is a {\it graph\/}; a 
simply
connected graph is a {\it simplicial tree}. A graph is a 
simplicial
tree if and only if it has no loops (topological 
embeddings of $S^1$\<).
Figure 3 gives the unique (up to isomorphism) trivalent 
simplicial tree.
One aspect of the negative curvature of trees is 
reflected in the fact
that the number of vertices of distance $\le n$ from a 
given vertex is
$3(2^n-1)+1$, a number which grows exponentially with $n$. 

The automorphism group of a simplicial tree $T$, $\a(T)$, 
is the group
of all self-homeomorphisms of $T$ which send vertices to 
vertices, edges
to edges, and which are linear on each edge. This is 
exactly the
automorphism group of the abstract complex 
$(V,E,\partial,\tau)$. An action
of a group $G$ on $T$ is a homomorphism from $G$ to 
$\a(T)$. The action
is said to be {\it without inversions\/} if for all $g\in 
G$ and all
$e\in E$ we have $e\cdot g\ne\ov  e$. To restrict to 
actions without

\fighere{9 pc}\caption{\smc Figure 2}

\fighere{13 pc}\caption{\smc Figure 3}

inversions is not a serious limitation since any action 
of $G$ on $T$
becomes an action without inversions on the first 
barycentric subdivision
of $T$ (obtained by splitting each edge of $T$ into its 
two halves and
adding a new vertex at the center of each edge). 

Given an action of a group $G$ on a simplicial tree $T$ 
there is the quotient
graph $\Gamma$. Its vertices are $V/G$; its oriented 
edges are $E/G$. The
geometric realization of this complex is naturally 
identified with the
quotient space $T/G$. 

\subheading{Graphs of groups} We are now ready to relate 
the theory of
group actions on simplicial trees to combinatorial group 
theory. The key 
to understanding the nature of an action of a group on a 
simplicial tree
is the notion of a graph of groups. A graph of groups is 
the following:
\roster
\item "(i)" A graph $\Gamma$, 
\item "(ii)" for each oriented edge or vertex $a$ of 
$\Gamma$ a group
$G_a$ such that if $e$ is an oriented edge then $G_{\ov 
e}=G_e$, and
finally,
\item "(iii)" if $v=i(e)$ then there is given an 
injective homomorphism
$G_e\hookrightarrow G_v$. 
\endroster
This graph of groups is said to be {\it over\/} $\Gamma$. 

There is the fundamental group of a graph of groups. A 
topological 
construction of the fundamental group goes as follows. 
For each edge or
vertex $a$ of $\Gamma$ choose a space $X_a$ whose 
fundamental group is
$G_a$. We can do this so that if $v$ is a vertex of $e$, 
then there is
an embedding $X_e\hookrightarrow X_v$ realizing the 
inclusion of groups.
We form the topological space $X(\Gamma)$ by beginning 
with the disjoint
union
$$
\coprod_{v\in V} X_v \coprod\coprod_{e\in E} X_e\times I,
$$
and (a) identifying $X_e\times I$ with $X_{\ov e}\times 
I$ via $(x,t)
\equiv (x,1-t)$ and (b) gluing $X_e\times\{0\}$ to 
$X_{i(e)}$ via the
given inclusion. The resulting topological space has 
fundamental group
which is the fundamental group of the graph of groups. A 
purely 
combinatorial construction of this fundamental group is 
given in 
\cite{20, pp.\ 41--42}. 

\ex{Example 1} Let $\Gamma$ be a single point. Then a 
graph over $\Gamma$
is simply a group. Its fundamental group is the group 
itself.
\endex

\ex{Example 2} Let $\Gamma$ be a graph with two vertices 
and a single
edge connecting them. Then a graph of groups over 
$\Gamma$ is the same
thing as an embedding of the edge group into two vertex 
groups. Its
fundamental group is the free product with amalgamation 
of the vertex
groups over the edge group. \endex

\ex{Example 3} Let $\Gamma$ be a graph with one vertex 
and a single 
edge, forming a loop. Then graph of groups over $\Gamma$ 
is a group $G_v$
and two embeddings $\varphi_0$ and $\varphi_1$ of another 
group $G_e$
into $G_v$. The fundamental group of this graph of groups 
is the
corresponding HNN-extension. A presentation of this 
extension is
$$
\langle G_v,s|s^{-1}\varphi_1(g) s=\varphi_0(g)\ 
\text{for all}\ 
g\in G_e\rangle\.
$$
\endex

 In general, the fundamental group of a graph of groups 
over a finite
graph can be described inductively by a finite sequence 
of operations
as in Examples 2 and 3. The fundamental group of an 
infinite graph of
groups is the inductive limit of the fundamental groups 
of the finite
subgraphs of groups.

Thus, it is clear that the operation of taking the 
fundamental group of
a graph of groups generalizes two basic operations of 
combinatorial
group theory---free product with amalgamation and 
HNN-extension.

There is a natural action of the fundamental group $G$ of 
a graph
of groups on a simplicial tree. To construct this action, 
let $X$ 
be the topological space, as described above, whose 
fundamental 
group is the fundamental group of the graph of groups. 
There is a 
closed subset $Y=\coprod_e X_e\times\{1/2\}\subset X$ 
which has a
collar neighborhood in $X$. Let $\wt X$ be the universal 
covering
of $X$, and let $\wt Y\subset\wt X$ be the preimage of 
$Y$. We define
the dual tree $T$ to $\wt Y\subset\wt X$. Its vertices 
are the
components of $\wt X-\wt Y$. Its unoriented edges are the 
components
of $\wt Y$. The vertices of an edge given by a component 
$\wt Y_0$
of $\wt Y$ are the two components of $\wt X-\wt Y$ which 
have $\wt Y_0$
in their closure. The fact the $\wt X$ is simply 
connected implies that
$T$ is contractible and hence is a simplicial tree.

There is a natural action of $G$ on $\wt X$. This action 
leaves $\wt Y$
invariant and hence defines an action of $G$ on $T$. It 
turns out that,
up to isomorphism, this action of $G$ on $T$ is 
independent of all the
choices involved in its construction. It is called {\it 
the universal
action associated with the graph of groups decomposition 
of $G$.} The
quotient $T/G$ is naturally identified with the original 
graph. Notice
that the stabilizer of a vertex or edge $a$ of $T$ is 
identified, up to
conjugation, with the group in the graph of groups 
indexed by the image
of $a$ in $T/G$. 

Here, we see for the first time the duality relationship 
between trees
and codimen\-sion-1 subsets.

\subheading{Structure theorem for groups acting on 
simplicial trees}
The main result in the theory of groups actions on 
simplicial trees is
a converse to this construction for a graph of groups 
decomposition of
$G$. It says:

\thm{Theorem 1} Let $G\times T\to T$ be an action without 
inversions of 
a group on a simplicial tree. Then there is a graph of 
groups over the
graph $T/G$ and an isomorphism from the fundamental group 
of this graph
of groups to $G$ in such a way that the action of $G$ on 
$T$ is identified
up to isomorphism with the universal action associated 
with the graph of
groups. \ethm

\thm{Corollary 2} Let $H$ act on a simplicial tree $T$. 
Suppose that no
point of $T$ is fixed by all $h\in H$. The $H$ has a 
HNN-decomposition
or has a nontrivial decomposition as a free product with
amalgamation.
\ethm

\demo{Proof} At the expense of subdividing $T$ we can 
suppose that $H$
acts without inversions. Consider the quotient graph 
$T/H$. We have a
graph of groups over $T/H$ whose fundamental group is 
identified with
$H$. There is a minimal subgraph $\Gamma\subset T/H$ such 
that the
fundamental group of the restricted graph of groups over 
$\Gamma$ 
includes isomorphically into $H$. This graph cannot be a 
single vertex
since the action of $H$ on $T$ does not fix any point. If 
$\Gamma$ has
a separating edge $e$, there is a nontrivial free product 
with 
amalgamation decomposition for $H$ as $H_a$ and $H_b$ 
amalgamated
along $G_e$, where $H_a$ and $H_b$ are the fundamental 
groups of the
graphs of groups over the two components of $\Gamma-e$. 
If $\Gamma$
has a nonseparating edge $e'$, then there is an 
HNN-decomposition for
$H$ given by the two embeddings of $G_e$ into the 
fundamental group
of the graph of groups over $\Gamma-e'$. \qed\enddemo

\ex{Example 4} Let $M$ be a manifold and let $N\subset M$ 
be a proper
submanifold of codimension 1. Let $\wt M$ be the 
universal covering of
$M$ and let $\wt N\subset \wt M$ be the preimage of $N$. 
Then $\wt N$
is collared in $\wt M$. We define a tree dual to $\wt 
N\subset\wt M$.
Its vertices are the components of $\wt M-\wt N$. Its 
edges are the
components of $\wt N$. We define the endpoints of an edge 
associated
to $\wt N_0$ to be the components of $\wt M-\wt N$ 
containing $\wt N_0$
in their closure. Since $\wt N$ is collared in $\wt M$, 
it follows that
each edge has two endpoints. Thus, we have defined a 
graph. Since $\wt M$\
is simply connected, this graph is a tree. The action of 
$\pi_1(M)$ on
$\wt M$ defines an action of $\pi_1(M)$ on $T$. The 
quotient is the
graph dual to $N\subset M$. \endex

\rem{Applications} As a first application notice that a 
group acts 
freely on a simplicial tree if and only if it is a free 
group. If $G$
acts freely on a tree, then $G$ is identified with the 
fundamental
group of the graph $T/G$, which is a free group. 
Conversely, if $G$ is
free on the set $W$, then the wedge of circles indexed by 
$W$ has
fundamental group identified with $G$. The group $G$ acts 
freely on
the universal covering of this wedge, which is a tree. 
\endrem

As our first application we have the famous

\thm\nofrills{Corollary 3\ \RM{(Schreier's theorem).}}\ 
Every subgroup of
a free group is free. \ethm

Using the relation between the number of generators of a 
free group and
the Euler characteristic of the wedge of circles, one can 
also establish
the Schreier index formula.

{\it If $G$ is a free group of rank $r$, and if 
$G'\subset G$ is a subgroup of
index $n$, then $G'$ is a free group of rank $n(r-1)+1$.}

One can also use the theory of groups acting on trees to 
give a generalization
of the Kurosh subgroup theorem. Here is the classical 
statement.

\thm\nofrills{Theorem \RM{(Kurosh subgroup theorem).}}\ 
Let $G=G_1*G_2$
be a free product. Let $G'\subset G$ be a subgroup. 
Suppose that $G'$
is not decomposable nontrivially as a free product. Then 
either $G'
\cong \bz$ or $G'$ is conjugate in $G$ to a subgroup of 
either $G_1$ or
$G_2$.\ethm

\demo{Proof} The decomposition $G=G_1*G_2$ gives an 
action of $G$ on a
simplicial tree $T$ with trivial edge stabilizers and 
with each vertex
stabilizer being conjugate in $G$ to either $G_1$ or 
$G_2$. Consider
the induced action of $G'$ on $T$. Then $G'$ has a graph 
of groups 
decomposition where all the vertex groups are conjugate 
in $G$ to
subgroups of either $G_1$ or $G_2$ and all the edge 
groups are trivial.
If this decomposition is trivial, then $G'$ is conjugate 
to a
subgroup of either $G_1$ or $G_2$. If it is nontrivial, 
then either 
$G'$ is a nontrivial free product or $G'$ has a 
nontrivial free factor.\qed\enddemo

The generalization which is a natural consequence of the 
theory of groups
acting on trees is:

\thm{Theorem 5} Let $G$ be the fundamental group of a 
graph of groups,
with the vertex groups being $G_v$ and the edge groups 
being $G_e$.
Suppose that $G'\subset G$ is a subgroup with the 
property that its
intersection with every conjugate in $G$ of each $G_e$ is 
trivial.
Then $G'$ is  the fundamental group of a graph of groups 
all of whose
edge groups are trivial. In particular, $G'$ is 
isomorphic to a free 
product of a free group and intersections of $G$ with 
various
conjugates of the $G_v$. \ethm

\demo{Proof} The restriction of the universal action of 
$G$ to $G'$
produces an action of $G'$ on a tree with trivial edge 
stabilizers
and with vertex stabilizers exactly the intersections of 
$G'$ with
the conjugates of the $G_v$.\qed\enddemo

\thm{Corollary 6} With $G$ and $G'$ as in Theorem \RM5, 
if $G'$ is
not a nontrivial free product, then it is either 
isomorphic to $\bz$
or is conjugate to a subgroup of $G_v$ for some vertex 
$v$. \ethm

\thm{Corollary 7} With $G$ as in Theorem \RM5, if 
$G'\subset G$ is a
subgroup with the property that the intersection of $G'$ 
with every
conjugate of $G_v$ in $G$ is trivial, then $G'$ is a free 
group. \ethm

\demo{Proof} The point stabilizers for the universal 
action of $G$ are
subgroups of conjugates of the $G_v$. Thus, under the 
hypothesis, the
restriction of the universal action to $G'$ is 
free.\qed\enddemo

These examples and applications give ample justification 
for the
statement that the theory of groups acting on trees is a 
natural
extension of classical combinatorial group theory.

\subheading{The tree associated to $SL_2$ over a local 
field}
One of the main reasons that Serre was led to investigate 
the
theory of groups acting on trees was to construct a space 
to play the 
role for $SL_2(\bq_p)$ that the upper half-plane plays 
for $SL_2(\br)$.
That space is a tree. Here is the outline of Serre's 
construction.

Let $K$ be a (commutative) field with a discrete 
valuation $v\:
 K^*\to\bz$.
Recall that this means that $v$ is a homomorphism from 
the multiplicative
group $K^*$ of the field onto the integers with the 
property that
$$
v(x+y)\ge\min(v(x),v(y))\.
$$
By convention we set $v(0)=+\infty$. The valuation ring 
$\scr O(v)$ is
the ring of $\{x\in K| v(x)\ge 0\}$. This ring is a local 
ring
with maximal ideal generated by any $\pi\in\scr O(v)$ 
with the property
that $v(\pi)=1$. The quotient $\scr O(v)/\pi\scr O(v)$ is 
called the
residue field $k_v$ of the valuation.

Let $W$ be the vector space $K^2$ over $K$. The group 
$SL_2(K)$ of
$2\times 2$\<-matrices with entries in $K$ and 
determinant 1 is
naturally the group of volume-preserving $K$\<-linear 
automorphisms
of $W$. 

An $\scr O(v)$\<-lattice in $W$ is a finitely generated 
$\scr O(v)$\<-submodule
$L\subset W$ which generates $W$ as a vector space over 
$K$. Such a module
is a free $\scr O(v)$\<-module of rank 2. We say that 
lattices $L_1$ and
$L_2$ are equivalent (homothetic) if there is $\alpha\in 
K^*$ such that
$L_1=\alpha\cdot L_2$. We define the set of vertices $V$ 
of a graph to be
the set of homothety classes of $\scr O(v)$\<-lattices of 
$W$. We join
two vertices corresponding to homothety classes having 
lattice
representatives $L_1$ and $L_2$ related in the following 
way: There is an $\scr O(v)$\<-basis
$\{e,f\}$ for $L_1$ so that $\{\pi e,f\}$ is an $\scr 
O(v)$\<-basis for
$L_2$. It is an easy exercise to show that this defines a 
simplicial tree 
on which $SL_2(K)$ acts. The stabilizer of any vertex is 
a conjugate
in $GL_2(K)$ of $SL_2(K)$. The quotient graph is the 
interval. A 
fundamental domain for the action on the tree is the 
interval connecting
the class of the lattice with standard basis $\{e,f\}$ to 
the class of
the lattice with basis $\{e,\pi f\}$. The stabilizer of 
the edge joining
these two lattices is the subgroup
$$
\Delta=\lf\{\pmatrix a & b\\ c & d\endpmatrix \in 
SL_2(\scr O(v))
| c\in \pi\scr O(v)\rt\}\.
$$
Thus, we have
$$
SL_2(K)\cong SL_2(\scr O(v))*_\Delta SL_2(\scr O(v))'
$$
where
$$
SL_2(\scr O(v))'=\pmatrix 1 & 0\\ 0 & \pi^{-1}\endpmatrix
SL_2(\scr O(v)) \pmatrix 1 & 0\\ 0 & \pi\endpmatrix\.
$$

If $A\in SL_2(K)$ and if $[L]\in T$, then the distance 
that $A$ moves
$[L]$ is given as follows. We take an $\scr O(v)$\<-basis 
for $L$ and
use it to express $A$ as a $2\times 2$ matrix. The 
absolute value of
the minimum of the valuation of the 4 matrix entries is 
the distance
that $[L]$ is moved. In particular, if trace $(A)$ has 
negative valuation
then $A$ fixes no point of $T$. 

\rem{Applications} In the special case of $SL_2(\bq_p)$ 
we have a
decomposition
$$
SL_2(\bq_p)=SL_2(\bz_p)*_\Delta SL_2(\bz_p)'\.
$$
The maximal compact subgroups of $SL_2(\bq_p)$ are the 
conjugates of
$SL_2(\bz_p)$. We say that a subgroup of $SL_2(\bq_p)$ is 
discrete if
its intersection with any maximal compact subgroup is 
finite. \endrem

\thm\nofrills{Corollary 8\ \RM{(Ihara's Theorem).}}\ Let 
$G\subset SL_2
(\bq_p)$ be a torsion-free, discrete subgroup. Then $G$ 
is a free group.
\ethm

\demo{Proof} If $G$ is both torsion-free and discrete, 
then its intersection
with any conjugate of $SL_2(\bz_p)$ is trivial. Applying 
Corollary 7 yields
the result.\qed\enddemo

If $C$ is a smooth curve, then let $\bc(C)$ denote the 
field of rational
functions on $C$. The valuations of $\bc (C)$ are in 
natural one-to-one
correspondence with the points of the completion 
$\widehat C$ of $C$.
For each $p\in \widehat C$ the associated valuation $v_p$ 
on $\bc(C)$ is
given by $v_p(f)$ is the order of the zero of $f$ at $p$. 
(If $f$ has a
pole at $p$, then by convention the order of the zero of 
$f$ at $p$ is
minus the order of the pole of $f$ at $p$.) Let $H$ be a 
fixed finitely
presented group. The representations of $H$ into 
$SL_2(\bc)$ form an 
affine complex algebraic variety, called the 
representation variety,
whose coordinate functions are the matrix entries of 
generators of $H$.
Suppose that we have an algebraic curve $C$ in this 
variety which is not
constant on the level of characters. That is to say there 
is $h\in H$
such that the trace of $\rho(h)$ varies as $\rho$ varies 
in $C$. Let $K$
be the function field of this curve. Each ideal point $p$ 
of $C$ is 
identified with a discrete valuation $v_p$ of $K$ 
supported at infinity 
in the sense that there is a regular (polynomial) 
function $f$ on $C$
with $v_p(f)<0$. Because of the condition that $C$ be a 
nontrivial 
curve of characters, for some ideal  point $p$, here is 
$h\in H$ such
that the regular function $\rho\mapsto 
\ot_h(\rho)=\te(\rho(h))$ has
negative value under $v_p$. 

Associated to such a valuation $v_p$ we have an action of 
$SL_2(K)$
on a simplicial tree $T$. We also have the tautological 
representation
of $H$ into $SL_2(K)$ (in fact into $SL_2$ of the 
coordinate ring of
regular functions on $C$\<). (In order to define this 
representation,
notice first that $C$ is a family of representations of 
$H$ into
$SL_2(\bc)$. The tautological representation assigns to 
$h\in H$
the matrix
$$
\pmatrix f_{11} & f_{12}\\ f_{21} & f_{22}\endpmatrix,
$$
where $f_{ij}(c)$ is the $ij$\<th entry of $c(h)\in 
SL_2(\bc)$.\<)
Thus, there is an induced action of $H$ on a simplicial 
tree.
If $v_p(\ot_h)<0$ it follows that the element $h$ fixes 
no point
of $T$. Thus, under our hypotheses on $C$ and $v_p$, it 
follows 
that the action of $H$ on the tree is nontrivial in the 
sense that
$H$ does not fix a point of the tree.

According to Corollary 2 we have proved

\thm\nofrills{Corollary 9 \RM{(\cite4).}}\ Let $H$ be a 
finitely presented
group. Suppose that the character variety of 
representations of $H$ into
$SL_2(\bc)$ is positive dimensional. Then there is an 
action of $H$ on a
tree without fixed point. In particular, $H$ has a 
nontrivial decomposition
as a free product with amalgamation or $H$ has an 
HNN-decomposition.\ethm

\heading 2. $\Lambda$\<-trees\endheading
In this section we define $\Lambda$\<-trees as the 
natural generalization
of simplicial trees. We prove the analogue of the result 
connecting trees
and $SL_2$. Namely, if $K$ is a local field with 
valuation $v\: K^*\to
\Lambda$ and if $G(K)$ is a rank-1 group, then there is a 
$\Lambda$\<-tree
on which $G(K)$ acts. We then take up the basics of the 
way groups act
on $\Lambda$\<-trees. We classify single isometries of 
$\Lambda$\<-trees
into three types and use this to define the hyperbolic 
length function of
an action. This leads to a definition of the space of all 
nontrivial,
minimal actions of a given group on $\Lambda$\<-trees. We 
finish the
section with a brief discussion of base change in 
$\Lambda$. 

Let us begin by reformulating the notion of a simplicial 
tree in a way
that will easily generalize. The vertices $V$ of a 
simplicial tree $T$
are a set with an integer-valued distance function. 
Namely, the distance
from $v_0$ to $v_1$ is the path distance or equivalently 
the minimal
number of edges in a simplicial path in $T$ from $v_0$ to 
$v_1$. The
entire tree $T$ can be reconstructed from the set $V$ and 
the distance
function. The reason is that two vertices of $T$ are 
joined by an edge
if and only if the distance between them is 1. A question 
arises as
to which integer-valued distance functions arise in this 
manner from
simplicial trees. The answer is not too hard to discover. 
It is based
on the notion of a segment. A $\bz$\<-segment in an 
integer valued
metric space is a subset which is isometric to a subset 
of the form
$\{t\in\bz| 0\le t\le n\}$. The integer $n$ is called the 
length
of the $\bz$\<-segment. The points corresponding to 0 and 
$n$ are 
called the endpoints of the $\bz$\<-segment.

\thm{Theorem 10} Let $(V,d)$ be an integer-valued metric 
space. 
Then there is a simplicial tree $T$  
such that the path distance function of $T$ is isometric 
to $(V,d)$ if and
only if the following hold
\roster
\item "(a)" For each $v,w\in V$ there is a 
$\bz$\<-segment in $V$
with endpoints $v$ and $w$. This simply means that there 
is a
sequence $v=v_0,v_1,\dots, v_n=w$ such that for all $i$ 
we have
$d(v_i,v_{i+1})=1$. 
\item "(b)" The intersection of two $\bz$\<-segments with 
an endpoint
in common is a $\bz$\<-segment.
\item "(c)" The union of two $\bz$\<-segments in $V$ 
whose  intersection
is a single point which is an endpoint of each is itself 
a $\bz$\<-segment.
\endroster
\ethm

Given a set with such an integer-valued distance 
function, one constructs
a graph by connecting all pairs of points at distance one 
from each other.
Condition (a) implies that the result is a connected 
graph. Condition (b)
implies that it has no loops. Condition (c) implies that 
the metric on $V$
agrees with the path metric on  this tree.

This definition can be generalized by replacing $\bz$ by 
any (totally) ordered
abelian group $\Lambda$. Before we do this, let us make a 
couple of 
introductory remarks about ordered abelian groups. An 
ordered abelian group
is an abelian group $\Lambda$ which is partitioned in 
three subsets $P$,
$N$, $\{0\}$ such that for each $x\ne 0$ we have exactly 
one of $x$,
$-x$ is contained in $P$ and with $P$ closed under 
addition. $P$ is said to be
the set of {\it positive\/} elements. We say that $x>y$ 
if $x-y\in P$.
A {\it convex subgroup\/} of $\Lambda$ is a subgroup 
$\Lambda_0$ with the
property that if $y\in\Lambda_0\cap P$ and if $0<x<y$ 
then $x\in\Lambda_0$.
The {\it rank\/} of an ordered abelian group $\Lambda$ is 
one less
than the length of the maximal chain of convex subgroups, 
each one proper
in the next. An ordered abelian group has rank 1 if and 
only if it is 
isomorphic to a subgroup of $\br$. 

A $\Lambda$\<-metric space is a pair $(X,d)$ where $X$ is 
a set and $d$
is a $\Lambda$\<-distance function $d\: X\times 
X\to\Lambda$ satisfying
the usual metric axioms. In $\Lambda$ there are segments 
$[a,b]$ given
by $\{\lambda\in\Lambda| a\le\lambda\le b\}$. More 
generally,
in a $\Lambda$\<-metric space a {\it segment\/} is a 
subset isometric
to some $[a,b]\subset\Lambda$. It is said to be {\it 
nondegenerate\/}
if $a<b$. As before, each nondegenerate segment has two 
endpoints.

\dfn\nofrills{Definition \RM{(\cite{24}, \cite{13}).}}\ A 
$\Lambda$\<-tree
is a $\Lambda$\<-metric space $(T,d)$ such that:
\roster
\item "(a)" For each $v,w\in T$ there is a 
$\Lambda$\<-segment in
$T$ with endpoints $v$ and $w$.
\item "(b)" The intersection of two $\Lambda$\<-segments 
in $T$ with
an endpoint in common is a $\Lambda$\<-segment. 
\item "(c)" The union of two $\Lambda$\<-segments of $T$ 
whose intersection
is a single point which is an endpoint of each is itself 
a $\Lambda$\<-segment.
\endroster
\enddfn

\ex{Example} Let $v\: K^*\to\Lambda$ be a possibly 
nondiscrete valuation.
Of course, by definition the value group $\Lambda$ is an 
ordered abelian
group. Associated to $SL_2(K)$ there is a 
$\Lambda$\<-tree on which it
acts. The points of this $\Lambda$\<-tree are again 
homothety classes
of $\scr O(v)$\<-lattices in $K^2$. The 
$\Lambda$\<-distance between
two homothety classes of lattices is defined as follows. 
Given the
classes there are representative lattices $L_0\subset 
L_1$ with quotient
$L_0/L_1$ and $\scr O(v)$\<-module of the form $\scr 
O(v)/\alpha\scr O(v)$
for some $\alpha\in\scr O(v)$. The distance between the 
classes is
$v(\alpha)$. In particular, the stabilizers of various 
points will be
conjugates of $SL_2(\scr O(v))$. Thus, every element 
$\gamma\in SL_2
(K)$ which fixes a point in this $\Lambda$\<-tree has 
$v(\ot(\gamma))\ge 0$.
If $\rho\: H\to SL_2(K)$ is a representation, then there 
is an induced action
of $H$ on the $\Lambda$\<-tree. This action will be 
without fixed point
for the whole group if there is $h\in H$ such that 
$v(\ot{\rho(h)})<0$. 
\endex

\subheading{The tree associated to a semisimple rank-1 
algebraic group
over a local field} The example in the previous section 
showed how
$SL_2$ over a local field with value group $\Lambda$ 
gives rise to a
$\Lambda$\<-tree. In fact this construction generalizes 
to any
semisimple real rank-1 group. We will content ourselves 
with considering
the case of $SO(n,1)$. Let $K$ be a field with a 
valuation $v\: K^*\to
\Lambda$. We suppose that the residue field $k_v$ is 
formally real in the
sense that $-1$ is not a sum of squares in $k_v$. (In 
particular, $k_v$
is of characteristic zero.) Let $q\: K^{n+1}\to K$ be the 
standard
quadratic form of type $(n,1)$; i.e.,
$$
q(x_0,\dots, x_n)=x_0x_1+x_2^2+\cdots+ x_n^2\.
$$
We denote by $\langle\cdot,\cdot\rangle$ the induced 
bilinear form.

Let $SO_K(n,1)\subset SL_{n+1}(K)$ be the automorphism 
group of $q$.
By a unimodular $\scr O(v)$\<-lattice in $K^{n+1}$ we 
mean a finitely
generated $\scr O(v)$\<-module which generates $K^{n+1}$ 
over $K$ and
for which there is a standard $\scr O(v)$\<-basis, i.e., 
an $\scr O
(v)$\<-basis $\{e_0,\dots, e_n\}$ with $q(\sum y_i 
e_i)=y_0y_1+y_2^2
+\cdots+ y_n^2$. Let $T$ be the set of unimodular $\scr 
O(v)$\<-lattices.
To define a $\Lambda$\<-metric on $T$ we need the 
following lemma which
is proved directly.

\thm{Lemma 11} If $L_0$ and $L_1$ are unimodular $\scr 
O(v)$\<-lattices,
then there is a standard $\scr O(v)$\<-basis for 
$L_0,\{e_0,e_1,\dots,
e_n\}$ and $\alpha\in \scr O(v)$ such that $$\{\alpha 
e_0,\alpha^{-1}
e_1,e_2,\dots, e_n\}$$ is a standard basis for $L_1$. 
\ethm

We then define the distance between $L_0$ and $L_1$ to be 
$v(\alpha)$.
It is not too hard to show that this space of unimodular 
$\scr O
(v)$\<-lattices with this $\Lambda$\<-distance function 
forms a
$\Lambda$\<-tree (see \cite{11}). Notice that the 
$\Lambda$\<-segment
between $L_0$ and $L_1$ is defined by taking the 
unimodular $\scr O
(v)$\<-lattices with bases $\{\beta e_0,\beta^{-1} 
e_1,e_2,\dots, e_n\}$
as $\beta$ ranges over elements of $\scr O(v)$ with 
$v(\beta)\le v(\alpha)$.
(This then generalizes the construction in Example 1 to 
$SO(n,1)$ and
we repeat to the other rank-1 groups.) 

The tree associated to $SO(n,1)$ over a local field $K$ 
will be a simplicial
tree exactly when the valuation on $K$ is discrete (i.e., 
when the value
group is $\bz$\<).

\subheading\nofrills{Basics of group actions on 
$\Lambda$\<-trees 
\RM{(cf.\ \cite{13, 1, 3}).}}\ We describe some of the 
basic results
about the way groups act on $\Lambda$\<-trees. The 
statements and proofs
are all elementary and directly ``tree-theoretic.'' Let 
us begin by
classifying a single automorphism $\alpha$ of a 
$\Lambda$\<-tree $T$.
There are three cases:
\roster
\item $\alpha$ has a fixed point in $T$. Then the fixed 
point set $F_\alpha$
of $\alpha$ is a subtree, and any point not in the fixed 
point set is
moved by $\alpha$ a distance equal to twice the distance 
to $F_\alpha$.
The midpoint of the segment joining $x$ to $\alpha(x)$ is 
contained in
$F_\alpha$. (See Figure 4.)
\item There is a segment of length $\lambda\in\Lambda$ 
but $\lambda\notin
2\Lambda$ which is flipped by $\alpha$. (See Figure 5.)
\item There is an axis for $\alpha$; that is to say there 
is an isometry
from a convex subgroup of $\Lambda$ into $T$ whose image 
is invariant under
$\alpha$ and on which $\alpha$ acts by translation by a 
positive amount
$\tau(\alpha)$.
\endroster

\fighere{4.5 pc}\caption{\smc Figure 4}
\fighere{3.5 pc}\caption{\smc Figure 5}
\fighere{6.5 pc}\caption{\smc Figure 6}

The set of all points in $T$ moved by $\alpha$ this 
distance $\tau(\alpha)$
themselves form the maximal axis $A_\alpha$ for $\alpha$. 
Any other point
$x\in T$ is moved by $\alpha$ a distance $2d(x,A_\alpha)+
\tau(\alpha)$
where $d(x,A_\alpha)$ is the distance from $x$ to 
$A_\alpha$. (See
Figure 6.)

The {\it hyperbolic length\/} $l(\alpha)$ of $\alpha$ is 
said to be 0
in the first two cases and $\tau(\alpha)$ in the last 
case. The
automorphism $\alpha$ is said to be {\it hyperbolic\/} if 
$\tau(\alpha)>0$,
{\it elliptic\/} if $\alpha$ has a fixed point and to be 
an {\it inversion\/}
in the remaining case. The {\it characteristic set\/} 
$A_\alpha$ of
$\alpha$ is $F_\alpha$ in Case 1, $A_\alpha$ in Case 3, 
and empty in 
Case 2. Thus, $C_\alpha$ is the set of points of $T$ 
which are moved
by $\alpha$ a distance equal to $l(\alpha)$. Clearly, 
there are no
inversions if $2\Lambda=\Lambda$. 

Notice that if $x\notin C_\alpha$, then the segment $S$ 
joining $x$
to $\alpha(x)$ has the property that $S\cap \alpha(S)$ is 
a segment
of positive length. We denote this by saying that the 
direction from
$\alpha(x)$ toward $x$ and the direction from $\alpha(x)$ 
toward
$\alpha^2(x)$ agree. Here is a lemma which indicates the 
`tree' 
nature of $\Lambda$\<-trees.

\thm{Lemma 12} Suppose that $g,h$ are isometries of a 
$\Lambda$\<-tree
$T$ with the property that $g,h$\RM, and $gh$ each have a 
fixed point in
$T$. Then there is a common fixed point for $g$ and 
$h$.\ethm

Notice that this result is false for the plane: two 
rotations of
opposite angle about distinct points in the plane fail to 
satisfy
this lemma.

\demo{Proof} Suppose that $g$ and $h$ have fixed points 
but that $F_g
\cap F_h=\varnothing$. Then there is a bridge between 
$F_g$ and $F_h$,
i.e., a segment which meets $F_g$ in one end and $F_h$ in 
the other.
Let $x$ be the initial point of this bridge. Then the 
direction from
$gh(x)$ toward $x$ and the direction from $gh(x)$ toward 
$(gh)^2(x)$\
are distinct. Thus, $x\in C_{gh}$. Since $x$ is not fixed 
by $gh$, it
follows that $gh$ is hyperbolic. (See Figure 7 on page 
100.)\qed\enddemo

\fighere{11.5 pc}\caption{\smc Figure 7}

\thm{Corollary 13} If $G$ is a finitely generated group 
of automorphisms 
of a $\Lambda$\<-tree such that each element in $G$ is 
elliptic, then
there is a point of $T$ fixed by the entire group $G$. In 
particular, if
$2\Lambda=\Lambda$ and if the hyperbolic length function 
of the action is
trivial, then there is a point of the tree fixed by the 
entire group.\ethm

Let us now describe some of the basic terminology in the 
theory of
$\Lambda$\<-trees. Let $T$ be a $\Lambda$\<-tree. A {\it 
direction\/}
from a point $x\in t$ is the germ of a nondegenerate 
$\Lambda$\<-segment
with one endpoint being $x$. The point $x$ is said to be 
a {\it branch
point\/} if there are at least three distinct directions 
from $x$. It
is said to be a {\it dead end\/} if there is only one 
direction. Otherwise,
$x$ is said to be a {\it regular point.}

If $T$ has a minimal action of a countable group $G$, 
then $T$ has no
dead ends and only countable many branch points. Each 
branch point
has at most countably many directions. It may well be the 
case when
$\Lambda=\br$ that the branch points are dense in $T$. 

If $G\times T\to T$ is an action of a group on a 
$\Lambda$\<-tree,
then the function associating to each $g\in G$, its 
hyperbolic
length, is a class function in the sense that it is 
constant on each
conjugacy class. We denote by $\scr C$ the set of 
conjugacy classes
in $G$. We define the {\it hyperbolic length function\/} 
of an
action of $G$ on a $\Lambda$\<-tree to be
$$
l\:\scr C\to\Lambda^{\ge 0}
$$
which assigns to each $c\in\scr C$ the hyperbolic length 
of any element
of $G$ in the class $c$. 

\subheading\nofrills{The space of actions 
\RM{(cf.\ \cite3).}}\ Several natural questions
arise. To what extent does the hyperbolic length function 
determine the
action? Which functions are hyperbolic length functions 
of actions?
The second question has been completely answered (at 
least for
subgroups of $\br$\<). There are some obvious necessary 
conditions
(each condition being an equation or weak inequality 
between the 
hyperbolic length of finitely many group elements), first 
laid out
in \cite3. In \cite{19} it was proved that these 
conditions characterize
the set of hyperbolic length functions.

Let us consider the first question. The idea is that the 
hyperbolic
length function of an action should be like the character 
of a 
representation. There are a couple of hurdles to surmount 
before this
analogy can be made precise. First of all, the inversions 
cause problems
much like they do in the simplicial case. Thus, as in the 
simplicial case,
one restricts to actions without inversions (which is no 
restriction at
all in the case $\Lambda=\br$\<). One operation which 
changes the action
but not its hyperbolic length function is to take a 
subtree. For most
actions there is a unique minimal invariant subtree. If 
the group is
finitely generated, the actions which do not necessarily 
have a unique
invariant subtree are those that fix points on the tree. 
These we call
{\it trivial actions.} It is natural to restrict to 
nontrivial actions
and to work with the minimal invariant subtree. In this 
context the
question then becomes to what extent the minimal 
invariant subtree of
a nontrivial action is determined by the hyperbolic 
length function 
of the action. There is a special case of 
`reducible-type' actions
much like the case of reducible representations where the 
hyperbolic
length function does not contain all the information 
about the minimal
invariant subtree, but this case is the exception. For 
all other
actions one can reconstruct the minimal action from the 
hyperbolic length
function. These considerations lead to a space of 
nontrivial, minimal
actions of a given finitely presented group on 
$\Lambda$\<-trees.
It is the space of hyperbolic length functions. When 
$\Lambda$ has a
topology, e.g., if $\Lambda\subset\br$ then this space of 
actions 
inherits a topology from the natural topology on the set 
of $\Lambda$\<-valued
functions on the group. For example, when $\Lambda=\br$ 
the space of
nontrivial, minimal actions is closed subspace of 
$(\br^{\ge 0})^{\scr C}
-\{0\}$. It is natural to divide this space by the action 
of $\br^+$ by
homotheties forming a projective space
$$
P((\br^{\ge 0})^{\scr C})=((\br^{ge0})^{\scr 
C}-\{0\}/\br^+)\.
$$
This projective space is compact, and the space of 
projective classes of
actions (or projectivized hyperbolic length functions) is 
a closed subset
of this projective space
$$
\scr PA(G)\subset P((\br^{\ge 0})^{\scr C}).
$$

\subheading\nofrills{Base change 
\RM{(cf.\ \cite1).}}\ Suppose that $\Lambda\subset\Lambda'$
is an inclusion of ordered abelian groups. Suppose that 
$T$ is a
$\Lambda$\<-tree. Then there is an extended tree 
$T\otimes_\Lambda\Lambda'$.
In brief one replaces each $\Lambda$\<-segment in $T$ 
with a
$\Lambda'$\<-segment. One example of this is 
$T\otimes_{\bz}\br$. This
operation takes a $\bz$ tree (which is really the set of 
vertices of a
simplicial tree) and replaces it with the $\br$\<-tree, 
which is the 
geometric realization of the simplicial tree. The 
operation $T
\otimes_{\bz} \bz[1/2]$ is the operation of barycentric 
subdivision.
Base change does not change the hyperbolic length function.

If $G$ acts on a $\Lambda$\<-tree $T$, then it acts 
without inversions
on the $\Lambda[1/2]$\<-tree 
$T\otimes_\Lambda\Lambda[1/2]$.

The operation of base change is a special case of a more 
general
construction that embeds $\Lambda$\<-metric spaces 
satisfying a
certain 4-point property isometrically into 
$\Lambda$\<-trees; see \cite1.

Finally, if $\Lambda'$ is a quotient of $\Lambda$ by a 
convex subgroup
$\Lambda_0$, then there is an analogous quotient 
operation that
applies to any $\Lambda$ tree to produce a 
$\Lambda'$\<-tree as quotient.
The fibers of the quotient map are $\Lambda_0$\<-trees.

\heading 3. Compactifying the space of 
characters\endheading
Let \<$G$\< be a finitely presented group. 
$R(G)\!=\!\hm(G,\!SO(n\!,\!1))$ is
naturally the real points of an affine algebraic variety 
defined
over $\bz$. In  fact, if $\{g_1,\dots, g_k\}$ are 
generators for 
$G$ then we have
$$
\hm(G,SO(n,1))\subset SO(n,1)^k \subset M(n\times 
n)^k=\br^{kn^2}
$$
is given by the polynomial equations which say that (1) 
all the $g_i$
are mapped to elements of $SO(n,1)$ and (2) that all the 
relations
among the $\{g_i\}$ which hold in $G$ hold for their 
images in $SO(n,1)$.
We denote by $SO_{\bc}(n,1)$ and $R_{\bc}(G)$ the complex 
versions of
these objects. Each $g\in G$ determines $n^2$ polynomial 
functions on
$R(G)$. These functions assign to each representation the 
matrix entries
of the representation on the given element $g$. The 
coordinate ring for
$R(G)$ is generated by these functions. These is an 
action of $SO_{\bc}
(n,1)$ and $R_{\bc}(G)$ by conjugation. The quotient 
affine algebraic
variety is called the character variety and is denoted 
$\chi_{\bc}(G)$
Though it is a complex variety, it is defined over $\br$. 
Its set of
real points, $\chi(G)$, is the equivalence classes of 
complex representations
with real characters. The polynomial functions on 
$\chi(G)$ are the
polynomial functions on $R(G)$ invariant by conjugation. 
In particular the
traces of the various elements of $G$ are polynomials on 
$\chi(G)$.
The character variety contains a subspace $Z(G)$ of 
equivalence classes
of real representations. The subspace $Z(G)$ is a 
semialgebraic subset
of $\chi(G)$ and is closed in the classical topology. 
Each fiber of
the map $R(G)\to Z(G)$ either is made up of 
representations whose images
are contained in parabolic subgroups of $SO(n,1)$ or is 
made up of
finitely many $SO(n,1)$\<-conjugacy classes of 
representations into
$SO(n,1)$.

The variety $\chi(G)$ and the subspace $Z(G)$ are usually 
not compact.
Our purpose here is produce natural compactifications of 
$Z(G)$ and
$\chi(G)$ as topological spaces and to interpret the 
ideal points
at infinity. The idea is to map the character variety 
into a projective
space whose homogeneous coordinates are indexed by $\scr 
C$, the set of
conjugacy classes in $G$. The map sends a representation 
to the point
whose homogeneous coordinates are the logs of the 
absolute values of its
traces on the conjugacy classes. The first result is that 
the image of
this map has compact closure in the projective space. As 
we go off to
infinity in the character variety at least one of these 
traces is going
to infinity. Thus, the point that we converge to in the 
projective space
measures the relative growth rates of the logs of the 
traces of the 
various conjugacy classes. It turns out that any such 
limit point in 
the projective space can be described by a valuation 
supported at infinity
on the character variety (or on some subvariety of it). 
These then
can be reinterpreted in terms of actions of $G$ on 
$\br$\<-trees.
This then is the statement for the case of $SO(n,1)$: A 
sequence of
representations of $G$ into $SO(n,1)$ which has unbounded 
characters
has a subsequence that converges to an action of $G$ on 
an $\br$\<-tree.
This material is explained in more detail in \cite{13} 
and \cite{11}.

\subheading{The projective space} We begin by describing 
the projective
space in which we shall work. Let $\scr C$ denote the set 
of conjugacy
classes in $G$. We denote by $P(\scr C)$ the projective 
space
$$
((\br^{\ge 0})^{\scr C}-\{0\})/\br^+
$$
where $\br^+$ acts by homotheties. Thus, a point in 
$P(\scr C)$ has 
homogeneous coordinates $[x_\gamma]_{\gamma\in \scr C}$ 
with the 
convention that each $x_\gamma$ is a nonnegative number.

\subheading{The main result} Before we broach all the 
technical details,
let us give a consequence which should serve to motivate 
the discussion.

\thm{Theorem 14} Let $\rho_k\: G\to SO(n,1)$ be a 
sequence of representations
with the property that for some $g\in G$ we have $\{\ot 
(\rho_k(g)\}_k$ 
is unbounded. Then after replacing the $\{\rho_k\}$ with 
a subsequence
we can find a nontrivial action
$$
\varphi\: G\times T\to T
$$
of $G$ on an $\br$\<-tree such that the positive part of 
the logs of the
absolute value of the traces of the $\rho_k(g)$ converge 
projectively
to the hyperbolic length function of $\varphi\!$\RM; 
i.e., so that if
$$
p_k=[\max(0,\log|\ot (\rho_k(\gamma))|)]_{\gamma\in\scr 
C} \in P(\scr C)
$$
then
$$
\lim_{k\mapsto \infty} 
p_k=[l(\varphi(\gamma))]_{\gamma\in\scr C}\.
$$
If all the representations $\rho_k$ of $G$ into $SO(n,1)$ 
are discrete
and faithful, then the limiting action of $G$ on an 
$\br$\<-tree has 
the property that for any nondegenerate segment $J\subset 
T$ the stabilizer
of $J$ under $\varphi$ is a virtually abelian group.
\ethm

The component of the identity $SO^+(n,1)$ of $SO(n,1)$ is 
an isometry group
of hyperbolic $n$\<-space. We define the hyperbolic 
length of $\alpha\in
SO(n,1)$ to be the minimum distance a point in hyperbolic 
space is
moved by $\alpha$. The quantity 
$\max(0,\log(|\ot(\alpha)|))$ differs
from the hyperbolic length of $\alpha$ by an amount 
bounded independent 
of $\alpha$. Expressed vaguely, the above result says 
that as actions of $G$
on hyperbolic $n$\<-space degenerate the hyperbolic 
lengths of these actions,
after rescaling, converge to the hyperbolic length of an 
action of $G$ on
an $\br$\<-tree.

The rest of this section is devoted to indicating how one 
establishes this
result.

\subheading{Mapping valuations into $P(\scr C)$} Let us 
describe how
valuations on the function field of $R(G)$ determine 
points of $P(\scr C)$.
Suppose that $\Lambda\subset\br$. Then any collection of 
elements 
$\{x_\gamma\}_{\gamma\in\scr C}$ of $\Lambda^{\ge 0}$, 
not all of which are
zero, determine an element of $P(\scr C)$. This can be 
generalized to any
ordered abelian group of finite rank. Suppose that 
$\Lambda$ is an ordered
abelian group and that $x,y\in\Lambda$ are nonnegative 
elements, not both
of which are zero. Then there is a well-defined ratio 
$x/y\in\br^{\ge 0}
\cup\infty$. More generally, if 
$\{x_\gamma\}_{\gamma\in\scr C}$ is a
collection of elements in $\Lambda$ with the property 
that $x_\gamma
\ge 0$ for all $\gamma\in\scr C$ and that $x_\gamma>0$ 
for some $\gamma
\in\scr C$ and if $\Lambda$ is of finite rank, then we 
have a well-defined
point
$$
[x_\gamma]_{\gamma\in\scr C}\in P(\scr C)\.
$$

Now suppose that $K$ is the function field of $R(G)$. 
This field contains
the trace functions $\ot_\gamma$ for all $\gamma\in\scr 
C$. A sequence
in $\chi(G)$ goes off to infinity if and only if at least
one of the traces $\ot_\gamma$ is unbounded on the 
sequence. 

\thm{Lemma 15} Let $v\: K^*\to \Lambda$ is a valuation, 
trivial on 
the constant functions, which is supported at infinity in 
the sense 
that for some $\gamma\in\scr C$ we have 
$v(\ot_\gamma)<0$. Then
$\Lambda$ is of finite rank. Thus we can define a point 
$\mu(v)$ in
$P(\scr C)$ by
$$
\mu(v)=[\max(0,-v(\ot_\gamma))]_{\gamma\in\scr C}\.
$$
Similarly, if $X\subset R(G)$ is a subvariety defined 
over $\bq$ and if
$v$ is a valuation on its function field $K_X$, trivial 
on the constant
functions, supported at infinity then we can define 
$\mu(v)$ to be the
same formula.\ethm

The reason that we take $\max(0,-v(\ot_\gamma))$ is that 
this number
is the logarithmic growth of $\ot_\gamma$ as measured by 
$v$. For
example, if the field is the function field of a curve 
and $v$ is a
valuation supported at an ideal point $p$ of the curve 
then $\max(0,
-v(f))$ is exactly the order of pole of $f$ at $p$. 

\subheading{How actions of $G$ on $\Lambda$\<-trees 
determine points in 
$P(\scr C)$} Another  source of points in $P(\scr C)$ is 
nontrivial
actions of $G$ on $\Lambda$\<-trees, once again for 
$\Lambda$ an ordered
abelian group of finite rank. Let $\varphi\: G\times T\to 
T$ be an 
action of $G$ on a $\Lambda$\<-tree. Let $T'$ be the 
$\Lambda[1/2]$\<-tree
$T\otimes_\Lambda \Lambda[1/2]$. Then there is an 
extended action of $G$
on $T'$ which is without inversions. Since $G$ is 
finitely generated, this
action is nontrivial if and only if its hyperbolic length 
function $l\:
\scr C\to\Lambda^{\ge 0}$ is nonzero. If the action is 
nontrivial, then
we have the point
$$
l(\varphi)=[l(\varphi(\gamma))]_{\gamma\in\scr C}\in 
P(\scr C)\.
$$

Actually, it is possible to realize the same point in 
$P(\scr C)$ by an
action of $G$ on an $\br$\<-tree. The reason is that we 
can find a
convex subgroup $\Lambda_0\subset\Lambda$ which contains 
all the
hyperbolic lengths and is minimal with respect to this 
property.
Let $\Lambda_1\subset\Lambda_0$ be the maximal proper 
convex subgroup
of $\Lambda_0$. Then the $T$ admits a $G$\<-invariant 
$\Lambda_0$\<-subtree.
This tree has a quotient $\Lambda_0/\Lambda_1$\<-tree on 
which $G$ acts.
Since $\Lambda_0/\Lambda_1$ is of rank 1, it embeds as a 
subgroup of $\br$.
Thus, by base change we can extend the $G$ action on this 
quotient tree
to a $G$ action on $\br$. This action determines the same 
point in 
$P(\scr C)$.

As we saw in the last section, the projective space $\scr 
PA(G)$ of
nontrivial minimal actions of $G$ on $\br$\<-trees sits 
by definition 
$P(\scr C)$. What we have done here is to show that any 
action of $G$
on a $\Lambda$\<-tree, for $\Lambda$ an ordered abelian 
group of finite 
rank, also determines a point of $\scr PA(G)\subset 
P(\scr C)$. 

There is a relationship between these two constructions 
of points in 
$P(\scr C)$, one from valuations on the function field of 
$\chi(G)$
and the other from actions on trees.

\thm{Theorem 16} Suppose that $X\subset R(G)$ is a 
subvariety defined over
$\bq$ and whose projection to $\chi(G)$ is unbounded. 
Suppose that $K_X$
is its function field and that $v\: K_X^*\to\Lambda$ is a 
valuation, trivial
on the constant functions, which is supported at infinity 
and which has
formally real residue field. Then associated to $v$ is an 
action of
$SO_K(n,1)$ on a $\Lambda$\<-tree. We have the 
tautological representation
of $G$ into $SO_K(n,1)$ and hence there is an induced 
action $\varphi$ of
$G$ on a $\Lambda$\<-tree. The image in $P(\scr C)$ of 
the valuation and
of this action of $G$ on the tree are the same\RM; i.e., 
$\mu(v)=l(\varphi)$.
\ethm

\subheading{Embedding the character variety in $P(\scr 
C)$} As we indicated
in the beginning of this section, the purpose for 
introducing the maps from
valuations and actions on trees to $P(\scr C)$ is to give 
representatives
for the ideal points of a compactification of the 
character variety,
$\chi(G)$. Toward this end we define a map 
$\theta\:\chi(G)\to P(G)$ by
$$
\theta([\rho])=[\max(0,\log(|\ot_\gamma(\rho)|))]_{\gamma%
\in\scr 
C}\.
$$
The map $\theta$ measures the relative sizes of the logs 
of the absolute
values of the traces of the images of the various 
conjugacy classes under
the representation. It is an elementary theorem (see 
\cite{13}) that the
image $\theta(\chi(G))$ has compact closure. Thus, there 
is an induced
compactification $\ov\chi(G)$. The set of ideal points of 
this
compactification, denoted by $B(\chi(G))$, is the compact 
subset of points
in $P(\scr C)$ which are limits of sequences 
$\{\theta([\rho_i])\}_i$ 
where $\{[\rho_i]\}_i$ is an unbounded sequence in 
$\chi(G)$. We denote
by $B(Z(G))$ the intersection of the closure of $Z(G)$ 
with $B(\chi(G))$.
Here are the two main results that were established in 
\cite{13} and
\cite{11}.

\thm{Theorem 17} For each point $b\in B(\chi(G))$ there 
exist a subvariety
$X$ of $R(G)$ defined over $\bq$ and a valuation $v$ on 
$K_X$, trivial
on the constant functions, supported at infinity such 
that $\mu(v)=b$. \ethm

\thm{Theorem 18} For each point $b\in B(Z(G))$ the 
valuation $v$
in Theorem \RM{17} with $\mu(v)=b$ can be chosen to have 
formally real
residue field. Thus, there is a nontrivial action 
$\varphi\: G\times T\to T$
on an $\br$\<-tree such that $b=l(\varphi)$. Said another 
way $B(Z(G))$ is
a subset of $\scr PA(G)\subset P(\scr C)$. Lastly, if the 
point $b\in B
(Z(G))$ is the limit of discrete and faithful 
representations, then the
action $\varphi$ has the property that the stabilizer of 
any nondegenerate
segment in the tree is a virtually abelian subgroup of 
$G$. \ethm

Theorem 14 is now an immediate consequence of this result.

There is an obvious corollary.

\thm{Corollary 19} Let $G$ be a finitely presented group 
nonvirtually
abelian group. If there is no nontrivial action of $G$ on 
an $\br$\<-tree
in which the stabilizer of each nondegenerate segment is 
virtually abelian,
then the space of conjugacy classes of discrete and 
faithful representations
of $G$ into $SO(n,1)$ is compact.\ethm

What is not clear in this result, however, is the meaning 
of the condition
about $G$ admitting no nontrivial actions on 
$\br$\<-trees with all
nondegenerate segments having virtually abelian 
stabilizers. This is a
question to which we shall return in \S5.

One thing is clear from the algebraic discussion, 
however. The discrete
valuations of a field are dense in the space of all 
valuations on the
field. Thus, it turns out that a dense subset of 
$B(\chi(G))$ is 
represented by discrete valuations. This leads to the 
following result.

\thm{Proposition 20} $B(Z(G))$ contains a countable dense 
subset represented
by nontrivial actions of $G$ on simplicial trees.
\ethm

From this result and Corollary 2 we have

\thm{Corollary 21} Let $G$ be a finitely generated group. 
If $Z(G)$ is not
compact\,\RM; i.e., if there is a sequence of 
representations $\{\rho_k\}_k$
of $G$ into $SO(n,1)$ such that for some $\gamma\in G$ 
the sequences of
traces $\{\ot(\rho_k(\gamma))\}_k$ is unbounded, then $G$ 
has either a
nontrivial free product with amalgamation decomposition 
or $G$ has an
HNN-decomposition. \ethm

\subheading{Geometric approach} What we have given here 
is an algebraic
approach to establish the degeneration of actions of $G$ 
on hyperbolic
space to actions of $G$ on $\br$\<-trees. There is, 
however, a purely
geometric approach to the same theory. In a geometric 
sense one can take
limits of hyperbolic space with a sequence of rescaled 
metrics so that
the curvature goes to $-\infty$. Any such limit is an 
$\br$\<-tree. 
Given a finitely generated group, there is a finite set 
of elements 
$\{g_1,\dots, g_t\}$ in $G$ such that for any $g\in G$ 
there is a 
polynomial $p_g$ in $t$ variables such that for any 
representation
$\rho\: G\to SO(n,1)$ the trace of $\rho(g)$ is bounded 
by the value
$p_g(\rho(g_1),\dots, \rho(g_t))$. Thus, given a sequence 
of representations
$\rho_k$ then for each $k$ we replace hyperbolic space by 
a constant 
$\lambda_k$ so that the maximum of the hyperbolic lengths 
of $\rho_k(g_1),
\dots, \rho_k(g_t)$ is 1. Any  limit of these rescaled 
hyperbolic spaces
will be an $\br$\<-tree on which there is an action of 
$G$. For more
details see \cite{12, \S\S 8--10}.

\heading4. Trees and codimensional-1 measured 
laminations\endheading

In this section we shall explain the relation between 
trees and 
codimension-1 transversely measured laminations. This 
relationship
generalizes Example 4 of \S1 where it was shown that if 
$N\subset M$
is a compact codimension-1 submanifold, then there is a 
dual action of 
$\pi_1(M)$ on a simplicial tree. 

\dfn{Definition} Let $M$ be a manifold. A (codimension-1 
transversely)
measured lamination $(\scr L,\mu)$ in $M$ consists of
\roster
\item "(a)" a closed subset $|\scr L|\subset M$ called the
{\it support of the lamination $\scr L$},
\item "(b)" a covering of $|\scr L|$ by open subsets $V$ 
of $M$,
called {\it flow boxes}, which have topological product 
structures
$V=U\times (a,b)$ (where $(a,b)$ denotes an open 
interval) so that
$|\scr L|\cap V$ is of the form $U\times X$ where 
$X\subset (a,b)$
is a closed subset, and
\item "(c)" for each open set as in (b) a Borel measure 
on the
interval $(a,b)$ with support equal to $X\subset(a,b)$.
\endroster
\enddfn

The open sets and measures are required to satisfy a 
compatibility
condition. We define the local leaves of the lamination 
in $V=U
\times (a,b)$ to be the slices $U\times\{x\}$ for $x\in 
X$. The 
first compatibility condition is that the germs of local 
leaves in
overlapping flow boxes agree. The second compatibility 
condition
involves transverse paths. A path in a flow box is said to 
{\it transverse to the lamination\/} if it is transverse 
to each
local leaf. Then the measures in the flow box can be 
integrated 
over transverse paths in that flow box to give a total 
measure.
If a transverse path lies in the intersection of two flow 
boxes 
then the total measures that are assigned to it in each 
flow 
box are required to agree. 

Two sets of flow boxes covering $|\scr L|$ define the 
same structure
if their union forms a compatible system of flow boxes. 
Equivalently,
we can view a measured lamination as a maximal family of 
compatible
flow boxes.

We define an equivalence relation on $|\scr L|$. This 
equivalence 
relation is generated by saying that two points are 
equivalent
if they both lie on the same local leaf in some flow box. 
The
equivalence classes are called the {\it leaves\/} of the 
measured
lamination. Each leaf is a connected codimension-1 
submanifold 
immersed in a 1-to-1 fashion in $M$. It meets each flow 
box in a
countable union of local leaves for that flow box. (See 
Figure 8.)

If the ambient manifold is compact, then the cross 
section of the support
of a measured lamination is the union of a cantor set 
where the measure
is diffuse, 

\fighere{9pc}\caption{\smc Figure 8}

\fighere{10 pc}\caption{\smc Figure 9}

\noindent an isolated set where the measure has 
$\delta$\<-masses and
nondegenerate intervals. If the lamination has no compact 
leaves, then 
every cross section meets the support in a cantor set.

The complement of the support of a measured lamination is 
an open subset
of $M$. It consists then of at most countably many 
components. It turns
out that when $M$ is compact, there can be countably many 
`thin' or
product components bounded by two parallel leaves. The 
other complementary
components are finite in number and are called the `big' 
complementary
regions. (See Figure 9.)
\ex{Example 1} A compact, codimension-1 submanifold 
$N\subset M$ is a
codimension-1 measured lamination; we assign the counting 
measure, that
is to say a $\delta$\<-mass of total mass 1 transverse to 
each component
of $N$.\endex

\ex{Example 2} Thurston \cite{23} has considered measured 
laminations
on a closed hyperbolic surface a genus $g$ all of whose 
leaves are 
geodesics. He has shown that the space of all such 
measured laminations
with the topology induced from the weak topology on 
measures is a real
vector space of dimension $6g-6$. The typical geodesic 
measured lamination
has $4g-4$ complementary regions each of which is an 
ideal triangle.
(See Figure 9.)\endex

It is always possible to thicken up any leaves that 
support $\delta$\<-masses
for the transverse measure to a parallel family of leaves 
with diffuse measure.
Let us assume that we have performed this operation, so 
that our measured
laminations have no $\delta$\<-masses.

Here is one theorem which relates actions on trees with 
measured 
laminations.

\thm\nofrills{Theorem 22 \RM{\cite{16}.}}\ Let $M$ be a 
compact manifold,
and let $(\scr L,\!\mu)$\< be a measured lamination in 
$M$. Suppose the
following hold\,\RM:
\roster
\item The leaves of the covering $\wt {\scr L}$ in the 
universal covering
$\wt M$ of $M$ are proper submanifolds.
\item If $x,y\in\wt M$, then there is a path joining them 
which is transverse
to $\wt{\scr L}$ and which meets each leaf of $\wt {\scr 
L}$ at most once.
\endroster
Then there is a dual $\br$\<-tree $T_{\scr L}$ to 
$\wt{\scr L}$ and an
action of $\pi_1(M)$ on $T_{\scr L}$. The branch points 
of $T_{\scr L}$
correspond to the big complementary regions of $\wt {\scr 
L}$ in $\wt M$,
and thus there are only finitely many branch points 
modulo the action of
$\pi_1(M)$. The regular points of $T_{\scr L}$ correspond 
to the leaves of
$\scr L$ which are not boundary components of ``big'' 
complementary
regions.

There is a continuous, $\pi_1(M)$ equivariant map $\wt 
M\to T_{\scr L}$
such that the inverse image of each branch point is a 
single complementary
region and the preimage of each regular point is either a 
single leaf
of the two leaves bounding a thin complementary 
component. \ethm

In this way we see that ``most'' measured laminations 
give rise to dual 
actions of the fundamental group on $\br$\<-trees. There 
is also a partial
converse to this result. Given an action of $\pi_1(M)$ on 
an $\br$\<-tree
$T$ it is possible to construct a transverse, 
$\pi_1(M)$\<-equivariant map 
from the universal covering $\wt M$ of $M$ to $T$. Using 
this map, one
can pull back a measured lamination from the tree and the 
metric on it. 
Thus, one produces a measured lamination associated to 
the action on the
tree. This measured lamination dominates the original 
action. For example,
if the lamination satisfies the hypothesis of Theorem 22, 
then dual to it
there is another action of $\pi_1(M)$ on an $\br$\<-tree. 
This new action
maps in an equivariant manner to the original action. If 
the original action
is minimal, then this map will be onto, and hence the new 
action dominates
the original one in a precise sense.

There are several difficulties with this construction. 
First, it is not
known that one can always do the construction so that the 
resulting 
lamination is dual to an action on a tree. Second, this 
construction
is not unique; there are many choices of the transverse 
map. In many
circumstances one hopes to find a ``best'' transverse 
map, but this
can be hard to achieve.

Let $M$ be a compact manifold. We say that an action of 
$\pi_1(M)$ on an
$\br$\<-tree is {\it geometric for\/} $M$ if there is a 
measured lamination
$(\scr L,\mu)$ in $M$ which satisfies the hypothesis of 
Theorem 22 and
such that the action of $\pi_1(M)$ dual to $(\scr L,\mu)$ 
is isomorphic
to the given action. We say that an action of an abstract 
finitely 
presented group $G$ is {\it geometric\/} if there is a 
compact manifold
$M$ with $\pi_1(M)=G$ and with the action geometric for 
$M$. One question
arises: {\it Is every minimal action of a finitely 
presented group
geometric\/}? If so it would follow that for a minimal 
action of a finitely
presented group $G$ on an $\br$\<-tree there are only 
finitely many branch
points modulo the action of $G$ and only finitely many 
directions from any
point $x$ of the tree modulo $G_x$. It is not known 
whether these statements
are true in general. 

In spite of these difficulties, measured laminations are 
an important
tool for studying actions of finitely presented groups on 
trees. One
extremely important result for measured laminations which 
gives a clue
as to the dynamics of actions on $\br$\<-trees is the 
following
decomposition result.

\thm\nofrills{Theorem 23\ \RM{(cf.\ \cite{13}).}}\ Let 
$(\scr L,\mu)$ be a
measured lamination in a compact manifold $M$. Then there 
is a decomposition
of $\scr L$ into finitely many disjoint sublaminations 
each of which has
support which is both open and closed in $|\scr L|$. Each 
of the sublaminations
is of one of the following three types\,\RM:
\roster
\item "(i)" a parallel family of compact leaves,
\item "(ii)" a twisted family of compact leaves with 
central member having
nontrivial normal bundle in $M$ and all the other leaves 
being two-sheeted
sections of this normal bundle, and 
\item "(iii)" a lamination in which every leaf is dense. 
\endroster
\ethm

We call a lamination of the last type an {\it exceptional 
minimal\/}
lamination.

One wonders if there is an analogous decomposition for 
actions of finitely
presented groups on $\br$\<-trees.

\heading 5. Applications and questions\endheading
In the analogy we have been developing between the 
theories of groups 
acting on simplicial trees and of groups acting on more 
general trees,
there is one missing ingredient. There is no analogue for 
the combinatorial
group theory in the simplicial case.There are some 
results that should be
viewed as partial steps in this direction. In this 
section, we shall discuss
some of these results and give applications of them to 
the spaces of hyperbolic
structures on groups. Two things emerge, clearly, from 
this discussion. First
of all, information about the combinatorial group 
theoretic consequences
of the existence of an action of $G$ on a 
$\Lambda$\<-tree has significant
consequences for the space of hyperbolic structures on 
$G$. Second,
most of the combinatorial group theoretic consequences to 
date have
followed from the use of measured laminations and related 
ergodic
considerations: first return maps, interval exchanges, etc.

\subheading{Applications to hyperbolic geometry}
Let us begin by formalizing the notion of a hyperbolic 
structure. Let
$G$ be a nonvirtually abelian, finitely generated group. 
For each $n$
we denote by $\scr H^n(G)$ the space of hyperbolic 
structures on $G$.
By definition this means the conjugacy classes of 
discrete and faithful
representations of $G$ into the isometry group of 
hyperbolic $n$\<-space.
Thought of another way a point in $\scr H^n(G)$ is a 
complete 
Riemannian $n$\<-manifold with all sectional curvatures 
equal to $-1$ and
with fundamental group identified with $G$.

\ex{Example 1} Let $G$ be the fundamental group of a 
closed surface of
genus at least 2. Then $\scr H^2(G)$ is the classical 
Teichm\"uller
space. It is homeomorphic to a Euclidean space of 
dimension $6g-6$.
\endex

\ex{Example 2} Let $G$ be the fundamental group of a 
closed hyperbolic
manifold of dimension $n\ge 3$. According to Mostow 
rigidity \cite{18},
$\scr H^n(G)$ consists of a single point. \endex

\ex{Example 3} Let $G$ be the fundamental group of a 
closed hyperbolic
manifold of dimension $n$. Then $\scr H^{n+1}(G)$ can be 
of positive dimension
(see \cite9).
\endex

The material in \S3 leads to the following result:

\thm\nofrills{Theorem 24 \RM{(\cite{13, 11}).}}\ There is 
a natural
compactification of
 $\scr H^n(G)$. Each ideal point of this compactification
is represented by a nontrivial action of $G$ on 
$\br$\<-trees with property
that the stabilizer of every nondegenerate segment of the 
tree is a 
virtually abelian subgroup of $G$. \ethm

In the special case of surfaces this gives a 
compactification of the
Tiechm\"uller space of a surface of genus $\ge 2$ by a 
space of actions
of the fundamental group on $\br$\<-trees. It turns out 
that all these
actions are geometric for the surface; i.e., the points 
at infinity in
this compactification can be viewed as geodesic measured 
laminations. For
more details see \cite{22} and \cite{13}.

The Teichm\"uller space of a surface is acted on properly 
discontinuously by 
the mapping class group (the group of outer automorphisms 
of the 
fundamental group of the surface). The action of the 
mapping class
group  extends to the compactification of Teichm\"uller. 
From this
one can deduce various cohomological results for the 
mapping class
group that make it look similar to an algebraic group. 
Motivated 
by this result, Culler-Vogtmann \cite5 have defined a 
space of free,
properly discontinuous actions of a free group on 
$\br$\<-trees. This
space is the analogue of the Teichm\"uller space for the 
outer 
automorphism group of the free group. The action of the 
outer
automorphism group of the free group on this space 
extends to the
compactification of this space of free, properly 
discontinuous
actions inside the projective space of all actions of the 
free 
group on $\br$\<-trees. By \cite3 the limit points are 
actions of
the free group on $\br$\<-trees with stabilizers of all 
nondegenerate
segments being virtually abelian. Using this action, 
Culler-Vogtmann
deduce many cohomological properties of the outer 
automorphism group
of a free group.

Theorem 24 leads immediately to the question: {\it Which 
groups act 
nontrivially on $\br$\<-trees with the stabilizer of 
every nondegenerate
segment in the tree being virtually abelian\/}?
The answer for simplicial trees follows immediately from 
the analysis in
\S1. 
\thm{Proposition 25} A group $G$ acts nontrivially on a 
simplicial tree
with all edge stabilizers being virtually abelian if and 
only if $G$
has either
\roster
\item "(A)" a nontrivial free production with 
amalgamation with the
amalgamating group being virtually abelian, or
\item "(B)" an HNN-decomposition with the subgroup being 
virtually
abelian. 
\endroster
\ethm

For fundamental groups of low dimensional manifolds, the 
answer is also
known. This is a consequence of a study of measured 
laminations in
3-manifolds and, in particular, the ability to do surgery 
on measured
laminations in 3-manifolds to make them incompressible.

\thm\nofrills{Theorem 26 \RM{(\cite{14, 15}).}}\ If $M$ 
is a
\RM3-manifold, then $\pi_1(M)$ acts nontrivially on an 
$\br$\<-tree
with the stabilizers of all nondegenerate segments being 
virtually
abelian if and only if it has such an action on a 
simplicial tree.\ethm

As an application we have

\thm{Corollary 27} Let $G$ be the fundamental group of a 
compact \RM3-manifold
$M$. Then for all $n$ the space $\scr H^n(G)$ is compact 
unless $G$ has a
decomposition of type \RM{(A)} or \RM{(B)} in Proposition 
\RM{25}.
 In particular,
if the interior of $M$ has a complete hyperbolic 
structure and has 
incompressible boundary, then $\scr H^n(G)$ is noncompact 
if and only if 
$M$ has an essential annulus which is not parallel into 
$\partial M$.
\ethm

One suspects that the first statement is true without the 
assumption that
$G$ is a 3-manifold group.

\rem{Conjecture} Let $G$ be a finitely presented group. 
Then for all
$n$ the space $\scr H^n(G)$ is compact unless $G$ has a 
decomposition
of type (A) or (B) in Proposition 25.\endrem

One also suspects that any action of $G$ on an 
$\br$\<-tree with 
segment stabilizers being virtually abelian can be 
approximated by
an action of $G$ on a  simplicial tree with the same 
property. This is
known for surface groups by \cite{22}, but not in general.

The fundamental group of a closed hyperbolic 
$k$\<-manifold, for $k\ge 3$,
has no decomposition of type (A) or (B). In line with the 
above suspicions
we have a result proved by the geometric rather than 
algebraic means.

\thm\nofrills{Theorem 28\ \RM{(\cite{12}).}}\ Suppose 
that $G$ is the
fundamental group of a closed hyperbolic manifold of 
dimension $k\ge 3$.
Then for all $n$, the space $\scr H^n(G)$ is compact.\ethm

\subheading{Combinatorial group theory for groups acting 
on $\br$\<-trees}
Perhaps the simplest question about the combinatorial 
analogues of the
simplicial case is the following: {\it which \RM(finitely 
presented\/\RM)
groups act freely on $\br$\<-trees\/}? Another question 
of importance,
suggested by the situation for valuations is: {\it Can 
one approximate
a nontrivial action of a group $G$ on an $\br$\<-tree by 
a nontrivial
action of $G$ on a simplicial tree\/\RM?
If the original action has stabilizers
of all nondegenerate segments being virtually abelian, it 
there a
simplicial approximation with this property\/\RM? \RM(For 
example, does the
boundary of Culler-Vogtmann space have a dense subset 
consisting of
simplicial actions\/\RM{?)}}

We finish by indicating some of the partial results 
concerning these two
questions. Because of the existence of geodesic 
laminations with simply
connected complementary regions, surface groups act 
freely on $\br$\<-trees
(see \cite{16}). Of course, any subgroup of $\br$ also 
acts freely on an
$\br$\<-tree. It follows easily that any free product of 
these groups acts
freely on an $\br$\<-tree. The question is whether there 
are any other groups
which act freely. By Theorem 26, for 3-manifold groups 
the answer is no.

We say that a group is {\it indecomposable\/} if it is 
not a nontrivial
free product. Clearly, it suffices to classify 
indecomposable groups which
act freely on $\br$\<-trees. In \cite{10} and \cite{17} 
the class of finitely
presented, indecomposable groups which have a nontrivial 
decomposition
as a free product with amalgamation or HNN-decomposition, 
where in each case
the subgroup is required to be virtually abelian, were 
studied. It was shown
that if such a group acts freely on an $\br$\<-tree, then 
it is either
a surface group or a free abelian group. This result is 
proved by invoking
the theory of measured laminations and results from 
ergodic theory. Another
result which follows from using the same techniques is: A 
minimal free action
of an indecomposable, noncyclic group $G$ on an 
$\br$\<-tree $T$ is mixing
in the following sense. Given any two nondegenerate 
segments $I$ and $J$ in 
$T$, there is a decomposition $I=I_1\cup\cdots\cup I_p$ 
and group elements
$g_i\in G$ such that for all $i$ we have $g_i\cdot 
I_i\subset J$. In 
particular, the orbit $G\cdot J$ of $J$ is the entire tree.

In a different direction, in \cite6 Gillet-Shalen proved 
that if $G$ acts
freely on an $\br$\<-tree which is induced by base change 
from a
$\Lambda$\<-tree where $\Lambda$ is a subgroup of $\br$ 
which generates
a rational vector space of rank at most 2, then $G$ is a 
free product of
surface groups and free abelian groups. By the same 
techniques they,
together with Skora \cite7, show that actions on such 
$\Lambda$\<-trees 
can be approximated by actions on simplicial trees 
(preserving the
virtually abelian segment stabilizer condition of 
relevant).

Recently, Rips has claimed that the only finitely 
generated groups which 
act freely on $\br$\<-trees are free products of surface 
groups and free
abelian groups.

\enddocument